\documentclass[12pt,fullpage, doublespace]{article}

\newif\ifpdf
\ifx\pdfoutput\undefined \pdffalse \else \pdfoutput=1 \pdftrue \fi

\ifpdf
\usepackage[pdftex]{graphicx}
\else
\usepackage{graphicx}

\fi

\usepackage{graphicx}
\usepackage{amsmath}
\usepackage{amsfonts}
\usepackage{epsfig}

\renewcommand{\baselinestretch} {1.3}

\makeatletter 
\def\singlespace{\def\baselinestretch{1}\@normalsize}

\title{{\sc \Large Moment properties of multivariate infinitely divisible laws and criteria for self-decomposability}}

        \author{{\em Theofanis ~Sapatinas} \footnote{Author for Correspondence: Email:
        \texttt{T.Sapatinas@ucy.ac.cy};
        Tel: \texttt{++357-22892653}; Fax: \texttt{++357-22892601.}}  \quad \text{and} \quad {\em Damodar N. ~Shanbhag},\\
         Department of Mathematics and Statistics,
         University of Cyprus,\\
         P.O. Box 20537,
         CY 1678 Nicosia,
         Cyprus.}
\date{}

\newcommand{\twofig}[4]
{ \hbox to\hsize{\hss
     \vbox{\psfig{figure=#1,width=#3,height=#4}}\qquad
     \vbox{\psfig{figure=#2,width=#3,height=#4}}
     \hss}
\vskip -0.0truein \hbox to\hsize{\hss
     \vbox{ \begin{center}\mbox{\footnotesize \hspace{0.0in} {(a)}
                      \hspace{#3} {(b)}  }  \end{center} }
     \hss}
\vskip 0.0truein }

\newcommand{\threefig}[5]{
\hbox to\hsize{\hss
     \vbox{\psfig{figure=#1,width=#4,height=#5}} \hspace{0.0in}
     \vbox{\psfig{figure=#2,width=#4,height=#5}} \hspace{0.2in}
     \vbox{\psfig{figure=#3,width=#4,height=#5}}
     \hss}
\vskip -0.0truein \hbox to\hsize{\hss
     \vbox{ \begin{center}\mbox{\footnotesize \hspace{0.0in} {(a)}
                      \hspace{#4} {(b)}   \hspace{#4} {(c)}}  \end{center} }
     \hss}
\vskip -0.1truein}

\newcommand{\fourfig}[6]
{ \hbox to\hsize{\hss
     \vbox{\psfig{figure=#1,width=#5,height=#6}}\qquad
     \vbox{\psfig{figure=#2,width=#5,height=#6}}
     \hss}
\vskip -0.0truein \hbox to\hsize{\hss
     \vbox{ \begin{center}\mbox{\footnotesize \hspace{0.1in} {(a)}
                      \hspace{#5} {(b)}  }  \end{center} }
     \hss}
\vskip 0.1truein \hbox to\hsize{\hss
     \vbox{\psfig{figure=#3,width=#5,height=#6}}\qquad
     \vbox{\psfig{figure=#4,width=#5,height=#6}}
     \hss}
\vskip -0.1truein
     \vbox{ \begin{center}\mbox{\footnotesize \hspace{0.1in} {(c)}
                      \hspace{#5} {(d)}  }  \end{center} }
\hbox to\hsize{\hss
     \hss}
\vskip -0.1truein }

\newcommand{\be}{\begin{equation}}
\newcommand{\ee}{\end{equation}}
\newcommand{\beqn}{\begin{eqnarray}}
\newcommand{\eeqn}{\end{eqnarray}}

\newcommand{\bt}{\beta}

\newcommand{\bft}{{\bf t}}

\newcommand{\PP}{\ensuremath{{\mathbb P}}}
\newcommand{\II}{\ensuremath{{\mathbb I}}}
\newcommand{\EE}{\ensuremath{{\mathbb E}}}

\newcommand{\RR}{\ensuremath{{\mathbb R}}}

\newcommand{\bbb}[1]{\boldsymbol #1}

\newcommand{\eee}{\mathrm{e}}
\newcommand{\iii}{\mathrm{i}}

\begin{document}

\ifpdf \DeclareGraphicsExtensions{.pdf, .jpg} \else
\DeclareGraphicsExtensions{.eps, .jpg} \fi

\maketitle


\begin{abstract}
Ramachandran (1969, Theorem 8) has shown that for any univariate
infinitely divisible distribution and any positive real number
$\alpha$, an absolute moment of order $\alpha$ relative to the
distribution exists (as a finite number) if and only if this is so
for a certain truncated version of the corresponding L$\acute{\rm
e}$vy measure. A generalized version of this result in the case of
multivariate infinitely divisible distributions, involving the
concept of g-moments, is given by Sato (1999, Theorem 25.3). We
extend Ramachandran's theorem to the multivariate case, keeping in
mind the immediate requirements under appropriate assumptions of
cumulant studies of the distributions referred to; the format of
Sato's theorem just referred to obviously varies from ours and seems
to be having a different agenda. Also, appealing to a further
criterion based on the L$\acute{\rm e}$vy measure, we identify in a
certain class of multivariate infinitely divisible distributions the
distributions that are self-decomposable; this throws new light on
structural aspects of certain multivariate distributions such as the
multivariate generalized hyperbolic distributions studied by
Barndorff-Nielsen (1977) and others. Various points of relevance to
the study are also addressed through specific examples.

\medskip

{\tt Keywords:} {\em Multivariate Generalized Hyperbolic
Distributions; Multivariate Indecomposability, Multivariate Infinite Divisibility; Multivariate
Self-Decomposability; Stable Distributions}

\medskip

{\em AMS (2000) Subject Classification:} {\sc  Primary 60E07;
Secondary 60E05; 60G51; 62H10}

\end{abstract}

\newpage

\section{\large Introduction}
\setcounter{equation}{0}

Infinite divisibility and their specialized versions, namely,
self-decomposability and stability, have generated considerable
interest among specialist in probability and statistics. There is
huge literature devoted to studies of these topics. Books such as
Lo$\grave{{\rm e}}$ve (1963), Linnik (1964), Feller (1966) and
Lukacs (1970) have been instrumental in providing the audience with
the basic material on these. More recent monographs such as
Bondesson (1992), Sato (1999) and Steutel \& van Harn (2004) have
unified and studied further contributions to the expanding
literature in this connection.

In view of Kendall \& Stuart (1963, Chapter 3), in which the
relations between moments and cumulants are addressed in detail, it
follows, under appropriate conditions, that the cumulants
corresponding to infinitely divisible distributions that exist have
some appealing features and have links with certain moments of
L$\acute{\rm e}$vy and Kolomogorov measures relative to these
distributions. Ramachandran (1969, Theorem 8) has shown, in the
univariate case, that for an infinitely divisible distribution the
existence of the absolute moment of order $\alpha \in (0, \infty)$
is equivalent to the existence of its analogue for a certain
truncated version of the corresponding L$\acute{\rm e}$vy measure.
(By a truncated version of a measure $\nu$ on $\RR^p$ ($p \geq 1$)
we mean the restriction of $\nu$ to some proper subset of $\RR^p$.)
This result plays a crucial role in studies related to cumulants of
infinitely divisible distributions, see, e.g., Gupta {\em et al.}
(1994) and Steutel \& van Harn (2004, Chapter IV, \S 7). Sato (1999,
Theorem 25.3) has given a multivariate extension of Theorem 8 of
Ramachandran (1969), involving the so-called $g$-moments. However,
it appears that Theorem 25.3 of Sato (1999) is not tailored to meet
the immediate needs for cumulant studies.

In a recent expository article, Gupta {\em et al.} (2009) have
unified the literature on infinitely divisible distributions with
special reference to moments and cumulants. In the process of doing
this, they have made several illuminating observations on the
behavior of cumulants of univariate and multivariate infinitely
divisible distributions, and have presented some new results in the
area. Gupta {\em et al.} (2009) also poses an open problem on a
multivariate extension of Theorem 8 of Ramachandran (1969). One of
the main tasks of the present article is to deal with this problem;
the problem that we have referred to here is of particular interest,
especially if one is concerned with aspects of cumulants of
multivariate infinitely divisible distributions. Interestingly, as a
by-product of our solution to the problem, it follows that for any
infinitely divisible distribution on $\RR^p$ ($p \geq 1$), under a
mild assumption, (in standard notation) the cumulant
$k_{r_1,\ldots,r_p}$ exists if the moment $\mu_{r_1,\ldots,r_p}$
exists (as a real number); in the univariate case, obviously, this
result holds without requiring the distribution to be infinitely
divisible.

In a somewhat different direction, there are questions relative to
structural aspects of the multivariate hyperbolic distributions of
Barndorff-Nielsen (1977) and their extensions with densities given
by equation (7.3) in the cited reference; each of the distributions
referred to here is indeed (in the notation of Barndorff-Nielsen
(1977)) a mixture of $N_n(\mu+u\beta\Delta, u\Delta)$ with respect
to $u$, where $u$ follows a certain generalized inverse Gaussian
distribution and $\mu$, $\beta$ and $\Delta$ are fixed with $\Delta$
nonsingular. (The extended versions have been termed the generalized
hyperbolic distributions, especially in the univariate case by,
e.g., Halgreen (1979, p. 14) and J{\o}rgensen (1982, p. 37).) As
claimed by Shanbhag \& Sreehari (1979, p. 24), there exist members
in the class of multivariate generalized hyperbolic distributions
that are not self-decomposable. Specific examples illustrating that
this is so can be found in, e.g., Pestana (1978, p. 54) and Rao \&
Shanbhag (2004, Example 3.3, Remark 3.4); Rao \& Shanbhag (2004)
consists of further information on the problem telling us, amongst
other things, that there exist members also in the smaller class of
multivariate hyperbolic distributions that are not
self-decomposable. In this article, we attempt a comprehensive
solution to a characterization problem that is linked with the
question on the structural aspects of multivariate hyperbolic and
multivariate generalized hyperbolic distributions being addressed.

The paper is organized as follows. In Section 2, a generalization of
Theorem 8 of Ramachandran (1969) to the case of multivariate
distributions is provided in conjunction with several relevant
observations and pertinent examples. In Section 3, a
characterization theorem, based on the property of
self-decomposability, is established for a certain class of mixtures
of multivariate distributions, and its implications are emphasized.
As mentioned before, the results of Section 2 are of importance in
cumulant studies and the results of Section 3 throw further light on
the structural aspects of multivariate generalized hyperbolic
distributions and some related distributions.

\section{\large Criteria based on L\'evy measure for the existence of moments for multivariate
infinitely divisible distributions} \setcounter{equation}{0}

From L$\acute{\rm e}$vy (1954), or any other appropriate source such
as Feller (1966, XVII.11) or Sato (1999, Theorem 8.1), it follows
that $\phi$ is the characteristic function (ch.f.) of an infinitely
divisible (i.d.) distribution on $\RR^p$ ($p \geq 1$) if and only if
it is of the form
\begin{equation}
\label{eq:id} \phi(\bbb t)= \exp \left\{ \iii <\bbb a, \bbb t>
-\frac{1}{2}Q(\bft) + \int_{\RR^p} \left(\eee^{\iii <\bbb t,\bbb x>}
-1 - \frac{\iii <\bbb t,\bbb x>}{1+\|\bbb x\|^2}\right)d\nu(\bbb x)
\right\}, \quad \bbb t \in \RR^p,
\end{equation}
where $\bbb a$ is a real vector, $Q$ is a nonnegative definite
quadratic form, and $\nu$ is a measure (referred to as L\'evy
measure) on the Borel $\sigma$-field of $\RR^p$ such that
$\nu(\{\bbb 0\})=0$ and
\begin{equation}
\label{eq:lm} \int_{\RR^p} \frac{\|\bbb x\|^2}{1+\|\bbb
x\|^2}\,d\nu(\bbb x) < \infty.
\end{equation}
It is easily seen that (\ref{eq:lm}) is equivalent to the condition
that
\begin{equation}
\label{eq:lme} \int_{\RR^p} \left(\min\{\|\bbb x\|,
\tau\}\right)^2\,d\nu(\bbb x) < \infty \quad \text{for any fixed}
\quad \tau \in (0, \infty).
\end{equation}
(Here, $<\cdot,\cdot>$ and $\|\cdot\|$ denote respectively the usual
inner product and usual norm on $\RR^p$.)

\medskip

As observed by Gupta {\em et al.} (2009), using essentially the
approach of Lo$\grave{{\rm e}}$ve (1963, Complement 9, p. 332), with
$\tau \in (0, \infty)$, we can rewrite (\ref{eq:id}) as
\begin{eqnarray}
\phi(\bbb t) & = & \exp \left\{ \iii <\bbb b, \bbb t>
-\frac{1}{2}Q(\bbb t) + \int_{\RR^p} \left(\eee^{\iii <\bbb t,\bbb
x>} -1 - \iii
<\bbb t,\bbb x>\right)d\nu_2(\bbb x) \right\} \nonumber \\
& & \times \exp \left\{-\nu_1(\RR^p) + \int_{\RR^p}\eee^{\iii <\bbb
t,\bbb x>}\,d\nu_1(\bbb x) \right\}, \quad \bbb t \in \RR^p,
\label{eq:ide}
\end{eqnarray}
where $\nu_1(\cdot) = \nu\left( \{ \bbb x:~\|\bbb x\| \geq \tau \}
\cap \cdot \right)$, $\nu_2 = \nu - \nu_1$, and (in obvious
notation)
$$
b_r = a_r + \int_{\RR^p} \left( \frac{x_r \|\bbb x\|^2}{1 + \|\bbb
x\|^2} \right)d\nu_2(\bbb x) - \int_{\RR^p} \left( \frac{x_r}{1 +
\|\bbb x\|^2} \right)d\nu_1(\bbb x), \quad r=1,2,\ldots,p.
$$
(Note that (\ref{eq:lm}), or the equivalent condition
(\ref{eq:lme}), trivially implies that $b_r$'s are well defined as
well as that $\nu$ is a $\sigma$-finite measure and $\nu_1$ is a
finite measure.)

\medskip

We begin now by giving the following theorem that extends Theorem 8
of Ramachandran (1969) to the case of i.d. distributions on $\RR^p$.
The theorem is clearly in a format which makes it easily applicable
in cumulant studies, especially when the constants $\alpha_r$ and
$\beta_r$, $r=1,2,\ldots,p$, defined in it are integers; the
expressions for cumulants corresponding to an i.d. distribution on
$\RR^p$ are related, under appropriate assumptions, to those for the
moments of measure $\nu_1$ referred to in (\ref{eq:ide}), as
observed, e.g., in Gupta {\em et al.} (2009).

\medskip

{\bf Theorem 1.} Let $\bbb X = (X_1,X_2,\ldots,X_p)$, $p \geq 1$, be
a $p$-component i.d. random vector with ch.f. $\phi$ satisfying (\ref{eq:id})
(and hence also (\ref{eq:ide})), and let
$\beta_r \in [0,\infty)$, $r=1,2,\ldots,p$. Then,
\begin{equation}
\label{eq:mom1} \EE \left( \prod_{r=1}^{p} |X_r|^{\alpha_r} \right)
< \infty \quad \text{for all} \quad \alpha_r \in [0, \beta_r], \quad
r=1,2,\ldots,p,
\end{equation}
if and only if
\begin{equation}
\label{eq:mom2} \int_{\RR^p} \left( \prod_{r=1}^{p}
|x_r|^{\alpha_r}\right)d\nu_1(\bbb x) < \infty \quad \text{for all}
\quad \alpha_r \in [0, \beta_r], \quad r=1,2,\ldots,p,
\end{equation}
where $\bbb x = (x_1,x_2,\ldots,x_p)$ and $\nu_1$ is as in
(\ref{eq:ide}).

\medskip

{\bf Proof.} That (\ref{eq:mom1}) implies (\ref{eq:mom2}) is an
obvious corollary to the first assertion of Theorem 5.1 of Gupta {\em
et al.} (2009). Now, to prove that (\ref{eq:mom2}) implies
(\ref{eq:mom1}), we may proceed as follows. Since $\bbb X$ is an
i.d. random vector with ch.f. $\phi$, in view of (\ref{eq:ide}), we
can see that $\phi$ is of the form
\begin{equation}
\label{eq:che} \phi(\bbb{t}) = \phi_2(\bbb{t}) \exp \{-\lambda +
\lambda\, \phi_1(\bbb{t}) \}, \quad \bbb{t} \in \RR^p,
\end{equation}
where $\lambda = \nu_1(\RR^p)$,
\begin{equation*}
\phi_1(\bbb t) = \left\{
\begin{array}{ll}
\lambda^{-1} \int_{\RR^p} \eee^{\iii
<\bbb t,\bbb x>}d\nu_1(\bbb x), \quad \bbb t \in \RR^p, &\mbox{if} \quad \lambda > 0, \\
1, \quad \bbb t \in \RR^p, &\mbox{if} \quad \lambda = 0
\end{array} \right.
\end{equation*}
and
$$
\phi_2({\boldsymbol t}) = \exp \left\{ \iii <\bbb b, \bbb t>
-\frac{1}{2}Q(\bbb t) + \int_{\RR^p} \left(\eee^{\iii <\bbb t,\bbb
x>} -1 - \iii <\bbb t,\bbb x>\right)d\nu_2(\bbb x) \right\}, \quad
\bbb t \in \RR^p.
$$

Essentially, adopting the approach of Feller (1966, p. 534) or
others such as Sato (1999, Lemmas 25.6, 25.7) or Steutel \& van Harn
(2004, Chapter IV, \S 7) given in the univariate case, it can be
easily seen that the distribution corresponding to the ch.f.
$\phi_2$ has a (full) moment sequence and, hence, satisfies the
analogue of (\ref{eq:mom1}). Since this disposes of the case of
$\lambda=0$ trivially, assume then that $\lambda >0$, where
$\lambda$ is as in (\ref{eq:che}). From (\ref{eq:che}), we see, in
this case, that
\begin{equation}
\label{eq:che2} \phi(\bbb t) = \phi_2(\bbb t) \sum_{j=0}^{\infty}
\eee^{-\lambda} \frac{\left( \lambda \phi_1(\bbb t)\right)^j}{j!},
\quad \bbb t \in \RR^p,
\end{equation}
which, in turn, implies (on noting, in particular, that $\phi$ is a
mixture of ch.f.'s $\phi_2(\bbb t)(\phi_1(\bbb t))^j$, $\bbb t \in
\RR^p$, for $j \in \{0,1,\ldots\}$, where the mixing distribution is
Poisson with mean $\lambda$) in view of Fubini's theorem that, for
each $\alpha_r \in [0, \beta_r]$, $r=1,2,\ldots,p$,
\begin{equation}
\label{eq:mom} \EE \left( \prod_{r=1}^{p} |X_r|^{\alpha_r} \right) =
\sum_{j=0}^{\infty} \eee^{-\lambda} \frac{\lambda^j}{j!}\, \EE
\left( \prod_{r=1}^{p} \left|X_{2r} + \sum_{k=1}^{j} X_{1r}^{(k)}
\right|^{\alpha_r} \right),
\end{equation}
where $(X_{21},X_{22},\ldots,X_{2p})$ and
$(X_{11}^{(k)},X_{12}^{(k)},\ldots,X_{1p}^{(k)})$, $k=1,2,\ldots,j$,
are mutually independent random vectors with the first one having
ch.f. $\phi_2$ and each of the remaining ones having ch.f. $\phi_1$.
(In (\ref{eq:mom}) we adopt the convention that the summations with
respect to $k$ equal 0 if $j=0$.)

Assume now that (\ref{eq:mom2}) holds. In view of (\ref{eq:mom}) and
(\ref{eq:mom2}), it easily follows that, for each $\alpha_r \in [0,
\beta_r]$, $r=1,2,\ldots,p$,
\begin{eqnarray}
\EE \left( \prod_{r=1}^{p} |X_r|^{\alpha_r} \right) & \leq &
\sum_{j=0}^{\infty} \eee^{-\lambda} \frac{\lambda^j}{j!}
(j+1)^{(\sum_{r=1}^{p}\alpha_r)}\, \EE \left\{ \prod_{r=1}^{p}
\sum_{k=0}^{j} \left| X_{1r}^{(k)}
\right|^{\alpha_r} \right\} \nonumber \\
& \leq & \sum_{j=0}^{\infty} \eee^{-\lambda} \frac{\lambda^j}{j!}
(j+1)^{(\sum_{r=1}^{p}\alpha_r)}\, \EE \left\{
\prod_{k=0}^{j}\prod_{r=1}^{p} \left( 1 + \left| X_{1r}^{(k)}
\right|^{\alpha_r} \right) \right\} \nonumber \\
& \leq & c \sum_{j=0}^{\infty} \eee^{-\lambda} \frac{(c
\lambda)^j}{j!} (j+1)^{(\sum_{r=1}^{p}\alpha_r)} < \infty,
\label{eq:bound1}
\end{eqnarray}
where
$$
c = \max \left\{  \EE \left( \prod_{r=1}^{p} \left(1+
\left|X_{1r}^{(1)}\right|^{\alpha_r}\right) \right),\EE \left(
\prod_{r=1}^{p} \left(1+\left|X_{2r}\right|^{\alpha_r}\right)
\right) \right\} < \infty,
$$
and, for notational convenience, we denote $X_{2r}$ by
$X_{1r}^{(0)}$ for each $r$. (We appeal here to, amongst other
things, the relevant independence of the vectors concerned and the
elementary inequality, for $\alpha \in [0, \infty)$, $y_k \in \RR$,
$k=0,1,\ldots,j$, that $|\sum_{k=0}^{j}y_k|^{\alpha} \leq
(j+1)^{\alpha} \max_{0 \leq k \leq j} |y_k|^{\alpha} \leq
(j+1)^{\alpha} \sum_{k=0}^{j} |y_k|^{\alpha}$, and also view the
last summation with respect to $j$ in (\ref{eq:bound1}) as the
expected value of a function of a Poisson random variable.)
Obviously, (\ref{eq:bound1}) yields the validity of (\ref{eq:mom1}).
This shows that (\ref{eq:mom2}) implies (\ref{eq:mom1}). Hence,
Theorem 1 follows. \hfill $\Box$

\medskip

The following corollary of Theorem 1 is essentially a version of
Theorem 8 of Ramachandran (1969); the cited author considers $\nu_1$
with ``$|x|>1$'' in place of ``$|x| \geq \tau$''.

\medskip

{\bf Corollary 1.} Let $X$ be an i.d. random variable with ch.f.
$\phi$, and let $\beta \in [0,\infty)$. Then,
\begin{equation*}
\label{eq:momr1} \EE \left( |X|^{\beta} \right) < \infty
\end{equation*}
if and only if
\begin{equation*}
\label{eq:momr2} \int_{\RR} |x|^{\beta}\,d\nu_1(x) < \infty,
\end{equation*}
where $\phi$ and $\nu_1$ are as in (\ref{eq:ide}) corresponding to
$p=1$.

\medskip

{\bf Proof.} The corollary follows readily from Theorem 1 by taking
$p=1$ because if $\mu$ is a finite measure on the Borel
$\sigma$-field of $\RR$ then $\int_{\RR}|x|^{\beta}\,d\mu(x) <
\infty$ if and only if $\int_{\RR} |x|^{\alpha}\,d\mu(x) < \infty$
for all $\alpha \in [0, \beta]$.  \hfill $\Box$

\medskip

{\bf Remark 1.} We point out that Theorem 1 is not covered by
Theorem 25.3 of Sato (1999) on $g$-moments (relative to the case of
$\tau=1$); that this statement is true is obvious from the example (i.e., Example 4)
given in Remark 5.5 of Gupta {\em et al.} (2009).  However, to
illustrate that neither the function $\prod_{r=1}^{p}
|x_r|^{\alpha_r}$, met in the statement of Theorem 1, nor its
modified version given by the function $\prod_{r=1}^{p}
|x_r|^{\alpha_r} \II(A)$, ${\bbb x} \in \RR^p$, $\tau
>0$, with $A= \{{\bbb x}:~\|{\bbb x}\| \geq \tau\}$ and $\II(A)$ as its indicator
function, is assured to be submultiplicative (where the terminology
refers to that of Definition 25.2 of Sato (1999, p. 159)), we may
consider the following example:

\medskip

{\bf Example 1.} Let $p=2$, $\bbb x=(x_1,x_2) \in \RR^2$ and $\tau
>0$, and consider
$$
g_1(\bbb x) = |x_1x_2|\II(A) \quad \mbox{and} \quad g_2(\bbb x) =
|x_1x_2|, \quad \bbb x \in \RR^2.
$$
For $\bbb x, \bbb y \in \RR^2$ such that $x_1 \geq \tau$, $y_2 \geq
\tau$, $x_2 = x_1^{-1}$ and $y_1 = y_2^{-1}$, we have
$$
g_r(\bbb x + \bbb y) = |(x_1+y_1)(x_2+y_2)| > x_1y_2, \quad r=1,2,
$$
and
$$
g_r(\bbb x) = g_r(\bbb y)=1, \quad r=1,2.
$$
Consequently, it is impossible to have here for each $r \in \{1,2\}$
a constant $a_r>0$ such that
$$
g_r(\bbb x+\bbb y) \leq a_r \,g_r(\bbb x) g_r(\bbb y), \quad r=1,2,
$$
for all $\bbb x, \bbb y \in \RR^2$. This supports the claim that we
have made above. (It is also now clear that the functions $\max \{
|x_1x_2|, c\}$ with $c>0$ are not submultiplicative.)

\medskip

{\bf Remark 2.} We may modify the example given in Remark 5.5 of
Gupta {\em et al.} (2009) to shed
further light on aspects of Theorem 1.  In particular, we can show,
with appropriate modifications to the example referred to, that if one or more
of certain $\beta_r$'s are positive then (\ref{eq:mom2}) with ``$[0,
\beta_{r}]$ (in respective places)'' replaced by ``$(0,
\beta_{r}]$'' does not imply (\ref{eq:mom1}) with the same change,
and demonstrate some curious phenomena of $(X_1,X_2)$ in this
connection. That this is so, is evident from the information
supplied by the following two examples:

\medskip

{\bf Example 2.} Let $(X_1,X_2)$ be an i.d. random vector with ch.f.
$\phi$ such that
$$
\phi(\bbb t) = \exp \{-\lambda + \lambda \psi(\bbb t) \}, \quad \bbb
t = (t_1, t_2) \in \RR^2,
$$
where $\lambda \in (0, \infty)$ and $\psi$ is the ch.f. of a random
vector $(Y_1,Y_2)$ satisfying
$$
(Y_1,Y_2) \stackrel{d}{=}  (V,V^{-1}),
$$
where $V$ is the modulus of a standard Cauchy random variable (and
``$\stackrel{d}{=}$'' denotes the equality in distribution).
Clearly, we have then for $\alpha_1, \alpha_2 \in [0,1]$ such that
$0 \leq |\alpha_1-\alpha_2|<1$, i.e., whenever $\alpha_1$, $\alpha_2$ both
lie in $(0, 1]$ (or both lie in $[0,1)$),
\begin{equation*}
\int_{\RR^2} |x_1|^{\alpha_1} |x_2|^{\alpha_2} \,d\nu_1(\bbb x) <
\int_{(0,\infty)^2} x_1^{\alpha_1} x_2^{\alpha_2} \,d\nu(\bbb x) =
\lambda \,\EE \left( V^{|\alpha_1-\alpha_2|} \right) < \infty.
\end{equation*}
However, in this case, we have
$$
\EE \left( |X_1X_2| \right) = \EE \left( X_1X_2 \right) = \infty
\quad \text{and} \quad \EE \left( |X_r| \right) = \EE \left( X_r
\right) = \infty, \quad r=1,2,
$$
violating the condition that $\EE \left( |X_1|^{\alpha_1}
|X_2|^{\alpha_2} \right)< \infty$ for all $\alpha_1, \alpha_2 \in
(0,1]$. (Appealing to Theorem 1, we can also see in this case that $
\EE \left( |X_1|^{\alpha_1} |X_2|^{\alpha_2} \right) < \infty$ for
each $\alpha_1, \alpha_2 \in [0,1)$.)

\medskip

{\bf Example 3.} Let $(X_1,X_2)$ be as in Example 2 with the
exception that $\psi$ in this case refers to the ch.f. of a random
vector $(Y_1,Y_2)$ satisfying
$$
(Y_1,Y_2) \stackrel{d}{=}
(V_{\gamma}^{\gamma},V_{\gamma}^{-\delta}),
$$
with $\gamma \in (0,1)$, $\delta \in (0, \infty)$ and $V_{\gamma}$
as a positive stable random variable with left extremity (i.e.,
$\inf \{x:~F_{V_{\gamma}}(x) > 0\}$, where $F_{V_{\gamma}}$ denotes
the distribution function (d.f.) of $V_{\gamma}$) zero and
characteristic exponent $\gamma$. It now follows that, for
$\alpha_1, \alpha_2 \in [0,1]$ with $\alpha_2 \neq 0$,
\begin{equation*}
\int_{\RR^2} |x_1|^{\alpha_1} |x_2|^{\alpha_2} \,d\nu_1(\bbb x)  <
\int_{(0,\infty)^2} x_1^{\alpha_1} x_2^{\alpha_2}\,d\nu(\bbb x) =
\lambda \, \EE \left( V_{\gamma}^{(\alpha_1 \gamma-\alpha_2 \delta)}
\right) < \infty.
\end{equation*}
(That the expectation appearing above is finite is seen, e.g., from  Bondesson (1992, p. 85) or
from Steutel \& van Harn (2004, p. 246).) Also, it is clear
that, in this case, $X_1$ and $X_2$ are nonnegative random variables
such that $\EE(X_1) = \EE(X_1X_2) = \infty$, with $X_2$ possessing a
(full) moment sequence.

\medskip

The following theorem enables us to understand the mechanism of
Theorem 1 better; the theorem addresses, amongst other things, the
problem posed in Remark 5.6 of Gupta {\em et al.} (2009), and its
proof that we have produced here is adapted partially from the proof
of Theorem 5.1 of Gupta {\em et al.} (2009).

\medskip

{\bf Theorem 2.} Let $\bbb X =(X_1,X_2,\ldots,X_p)$ and $\beta_r$,
$r=1,2,\ldots,p$, be as in Theorem 1, but for that there is an additional restriction now
that $\PP(X_r=0)<1$, $r=1,2,\ldots,p$. Then, (\ref{eq:mom1}) is
equivalent to the condition that
\begin{equation}
\label{eq:mom1Th2} \EE \left( \prod_{r=1}^{p} |X_r|^{\beta_r}
\right) < \infty.
\end{equation}
Moreover, we now have
\begin{equation}
\label{eq:propTh2} \PP \left\{ \prod_{r=1}^{p} |X_r| > 0 \right\} >
0.
\end{equation}

\medskip

{\bf Proof.} Trivially, the theorem is true for $p=1$. We shall
follow the method of induction with respect to $p$ to prove it,
noting that each subvector of $\bbb X$ is i.d. Assume then that $p >
1$ and the theorem holds in the case when $p'$, with $p' \in
\{1,2,...,p-1\}$, appears in place of $p$ in (\ref{eq:mom1Th2}) and
(\ref{eq:propTh2}). The theorem is clearly valid if $\bbb X$ is
$p$-variate normal and hence it is sufficient if we show that the
theorem holds when $\nu$ is non-null. Consequently, we can assume,
without loss of generality, that $\nu_1(\RR_{+}^p)>0$, where
$\RR_{+}=[0,\infty)$, and observe that the theorem follows if, under
the assumption, it is proved just that (\ref{eq:propTh2}) is valid
and, in view of Theorem 1 and the prevailing symmetry, that
(\ref{eq:mom1Th2}) implies (\ref{eq:mom2}) with $\RR^p$ replaced by
$\RR_{+}^p$.

Let $\lambda=\nu_1(\RR_{+}^p) >0$; then (\ref{eq:ide}) implies that,
for $n=1,2,\ldots$, we have $p$-variate ch.f.'s $\phi_2^{(n)}$ and
$\phi_3^{(n)}$ such that
\begin{equation}
\label{eq:chfth21} \phi(\bbb t) = \gamma_n\, \phi_1(\bbb t)
\phi_2^{(n)}(\bbb t) + (1-\gamma_n)\,\phi_3^{(n)}(\bbb t), \quad
\bbb t \in \RR^p,
\end{equation}
where
\begin{eqnarray*}
\phi_2^{(n)}(\bbb t) & = & \phi_2^{(1)}(\bbb t) \left(\phi_1(\bbb t)
\right)^{n-1}, \quad \bbb t \in \RR^p, \\
\phi_2^{(1)}(\bbb t) & = & \frac{\phi(\bbb t)}{\exp \{-\lambda +
\lambda\,\phi_1(\bbb t) \}}, \quad \bbb t \in \RR^p,\\
\phi_1(\bbb t) & = & \lambda^{-1}  \int_{\RR_{+}^p}  \eee^{\iii
<\bbb t, \bbb x>}\,d\nu_1(\bbb x), \quad \bbb t \in \RR^p,
\end{eqnarray*}
with $\gamma_n = \eee^{-\lambda} \lambda^n/n!$. 
(Note that $\lambda$ and $\phi_1$ considered here are not implied to be the same as those met in the proof of Theorem 1.)

Letting $(X_{11},X_{12},\ldots,X_{1p})$ and
$(X_{21}^{(n)},X_{22}^{(n)},\ldots,X_{2p}^{(n)})$, $n=1,2,\ldots$,
denote independent random vectors with ch.f.'s $\phi_1$ and
$\phi_2^{(n)}$, $n=1,2,\ldots$, respectively, we get from
(\ref{eq:chfth21}) that
\begin{equation}
\label{eq:mewmom1} \EE \left( \prod_{r=1}^{p} |X_r|^{\beta_r}
\right) \geq \gamma_n \, \EE \left( \prod_{r=1}^{p}
\left|X_{1r}+X_{2r}^{(n)}\right|^{\beta_r} \right), \quad
n=1,2,\ldots,
\end{equation}
and (noting especially that, for each $n$, $\phi_2^{(n+1)}(\bbb t) =
\phi_1(\bbb t) \phi_2^{(n)}(\bbb t)$, $t \in \RR^p$)
\begin{equation}
\label{eq:mewmom2} \PP \left\{ \prod_{r=1}^{p} |X_r| > 0 \right\}
\geq \gamma_n \, \PP\left\{ \prod_{r=1}^{p}
\left|X_{2r}^{(n+1)}\right|
> 0 \right\}, \quad n=1,2,\ldots\,.
\end{equation}
If $k<p$ of the $X_{1r}$, $r=1,2,\ldots,p$, are equal to zero almost
surely and the remainder satisfy the condition that $\PP\{X_{1r} >
0\} >0$, then we can take without loss of generality that
$\PP\{X_{11}=X_{12}=\ldots=X_{1k}=0\}=1$ and $\PP\{X_{1r} > 0\}
>0$, $r=k+1,k+2,\ldots,p$; we take here $k=0$ to mean that
$\PP\{X_{1r} >0\}>0$ for all $r \in \{1,2,\ldots,p\}$. For a
sufficiently large integer $n_0$, the distribution corresponding to
$\phi_2^{(n_0)}$ has at least one support point, say
$(c_1,c_2,\ldots,c_p)$ with $c_r \neq 0$ for each $r \in
\{1,2,\ldots,p\}$ and $c_r > 0$ for each $r \in
\{k+1,k+2,\ldots,p\}$; to see this note that, by assumption,
\begin{equation}
\label{eq:mewmom3} \PP \left\{ \prod_{r=1}^{k} |X_{2r}^{(1)}| > 0
\right\} > 0 \quad \text{if} \quad k \geq 1,
\end{equation}
and $\PP\{X_{11}=X_{12}=\ldots=X_{1k}=0\}=1$ together with
$\PP\{X_{1r} > 0\} >0$, $r=k+1,k+2,\ldots,p$. (Obviously,
(\ref{eq:mewmom3}) is a consequence of the assumption in the
inductive argument, especially because we have now that, for each
$n$, $(X_{21}^{(n)},X_{22}^{(n)},\ldots,X_{2k}^{(n)})
\stackrel{d}{=} (X_1,X_2,\ldots,X_k)$.)

In view of the aforementioned observation on the existence of the
support point $(c_1,c_2,\ldots,c_p)$ of the distribution relative to
$\phi_2^{(n_0)}$ with the stated properties, it follows that
\begin{equation}
\label{eq:fanisNewN} \PP \{ X_{2r}^{(n_0)} \in A_r, \,
r=1,2,\ldots,p \} >0,
\end{equation}
where
\begin{equation}
\label{eq:galisNewN} A_r = \left\{
\begin{array}{ll}
\left(\frac{3c_r}{2}, \frac{c_r}{2}\right) &\mbox{if} \quad c_r < 0, \\
\left(\frac{c_r}{2}, \frac{3c_r}{2}\right) &\mbox{if} \quad c_r > 0.
\end{array} \right.
\end{equation}
Consequently, by (\ref{eq:mewmom2}) we get that (\ref{eq:propTh2})
is valid, and, by (\ref{eq:mewmom1}) in conjunction with
(\ref{eq:mom1Th2}), that, for some constant $\eta \in (0,\infty)$
and $A_r$'s as in (\ref{eq:galisNewN}),
\begin{equation}
\label{eq:mewmom4} \infty > \EE \left( \prod_{r=1}^{p}
|X_r|^{\beta_r} \right) > \eta \, \EE \left( \prod_{r=1}^{p}
\left|X_{1r}+X_{2r}^{(n_0)}\right|^{\beta_r} \big|\, X_{2r}^{(n_0)}
\in A_r,\, r=1,2,\ldots,p \right).
\end{equation}
In view of the properties of $(X_{11},X_{12},\ldots,X_{1p})$ and
$(X_{21}^{(n_0)},X_{22}^{(n_0)},\ldots,X_{2p}^{(n_0)})$,
(\ref{eq:mewmom4}) implies then that
\begin{equation}
\label{eq:mewmom555} \EE \left( \prod_{r=1}^{p}
\left|X_{1r}+c_r^{\star}\right|^{\beta_r} \right) < \infty,
\end{equation}
where $c_r^{\star} = |c_r|/2$, $r=1,2,\ldots,p$. Since, for
$\alpha_r \in [0,\beta_r]$, $r=1,2,\ldots,p$,
$$
\left|X_{1r}\right|^{\alpha_r} \leq (c_r^{\star})^{\alpha_r-\beta_r}
\left|X_{1r}+c_r^{\star}\right|^{\beta_r}, \quad r=1,2,\ldots,p,
$$
it is hence obvious from (\ref{eq:mewmom555}) that
\begin{equation*}
\label{eq:mewmom6} \EE \left( \prod_{r=1}^{p} |X_{1r}|^{\alpha_r}
\right) < \infty \quad \text{for all} \quad \alpha_r \in [0,
\beta_r], \quad r=1,2,\ldots,p,
\end{equation*}
asserting that, as required for the completion of the proof of the
theorem, (\ref{eq:mom2}) with $\RR_{+}^p$ in place of $\RR^p$, is
met.  Hence, Theorem 2 follows. \hfill $\Box$

\medskip

{\bf Remark 3 (i).} Theorem 2 is not valid if the assumption of i.d. is
dropped. This is obvious from the following examples in which we use 
an indirect approach based on the first assertion of Theorem 2 to ascertain 
that the $\bbb X$ considered are non-i.d.:

\medskip

{\bf Example 4A.} Let $V$ and $V^{\star}$ be independent random variables
such that $\EE(|V|)=\infty$ and $V^{\star}$ is $\{0,1\}$-valued Bernoulli.
Define $X_1=VV^{\star}$ and $X_2=V(1-V^{\star})$. Clearly, the random vector $\bbb
X=(X_1,X_2)$ is such that $\PP\{ X_r=0 \}<1$, $r=1,2$, with $\PP
\{X_1X_2=0\}=1$ and hence with $\EE(|X_1X_2|)=0<\infty$. Also, in
this case, obviously $\EE(|X_1|) = \EE(|X_2|)= \infty$. That $\bbb
X$ is non-i.d. and the claim of Remark 3 is valid follows then
trivially from the first assertion of Theorem 2. However, it is
interesting to observe that in this example, we have
$\EE(|X_1|^{\alpha_1}|X_2|^{\alpha_2})=0 < \infty$ for all
$\alpha_1, \alpha_2 \in (0,1]$.

\medskip

{\bf Example 4B.} Let $\gamma \in (0,1)$ and $\bbb X=(X_1,X_2)$ be a
random vector such that $X_1$ is a positive random variable with
$\EE(X^\gamma)=\infty$ and $X_2=X_1^{-1}$ almost surely. Then,
clearly, we have $\EE(|X_1||X_2|)=\EE(X_1 X_2)=1 < \infty$. However,
in this case, it is not even true that
$\EE(|X_1|^{\alpha_1}|X_2|^{\alpha_2}) < \infty$ for all $\alpha_1,
\alpha_2 \in (0,1)$, since $\EE(|X_1|^{\gamma+\delta}|X_2|^{\delta})
= \EE(X_1^{\gamma})=\infty$ if $\delta \in (0, 1-\gamma)$. That
$\bbb X$ is non-i.d. and the claim of Remark 3 is valid follows
again trivially from the first assertion of Theorem 2.

\medskip

{\bf Remark 3 (ii).}  That the $\bbb X$ vectors dealt with in Examples 4A and 4B are non-i.d. can also be shown via 
alternative approaches without involving the findings of Theorem 2; in the remainder of this remark, and in Remark 3 (iii),
we illustrate as to why this is so. Let $(U,W)$ be a 2-component random vector with nonnegative and nondegenerate
components $U$ and $W$ such that $\PP\{UW=c\}=1$ for some (nonnegative) constant $c$. Since, in this case,
as a simple corollary to Theorem 2 of Shanbhag (1988), it follows that the distribution of $(U, (UW)^{1/2}, W)$ is indecomposable,
it is obvious then that the distribution of $(U,W)$ is indecomposable; to see this, note, in particular, that 
$(UW)^{1/2}$ is degenerate.  Consequently, we have the distribution of $\bbb X$ in Example 4A in the case  when $V$ 
is nonnegative and that of $\bbb X$ in Example 4B to be indecomposable and hence non-i.d. 

\medskip

{\bf Remark 3 (iii).} Let $(U,W)$ be as in Remark 3 (ii), but for a modification that $U$ and $W$ in this case are not 
necessarily nonnegative and also that $c$ is allowed here to be negative. Applying essentially a simpler version of the 
argument used in Shanbhag (1988) to prove its Theorem 2, one can see that the distribution of $(U,W)$ is decomposable if 
and only if for some $b_1, b_2 \neq 0$ with $b_1^2-4 c b_1 b_2^{-1} >0$, and some $\alpha, \beta \in (0,1)$,
\begin{equation*}
\label{eq:remark3(iii)} P\{ (U,W)= \bbb x\} = \left\{
\begin{array}{ll}
\alpha \beta &\mbox{if} \quad \bbb x=(a_1, a_2 b_1^{-1} b_2), \\
(1-\alpha) \beta &\mbox{if} \quad \bbb x=(a_2, a_1 b_1^{-1} b_2), \\
\alpha (1-\beta) &\mbox{if} \quad \bbb x=(-a_2, -a_1 b_1^{-1} b_2), \\
(1-\alpha) (1-\beta) &\mbox{if} \quad \bbb x=(-a_1, -a_2 b_1^{-1} b_2),
\end{array} \right.
\end{equation*}
where 
$$
a_1 = 2^{-1} \Big(b_1+\sqrt{b_1^2-4 c b_1 b_2^{-1}}\Big) \quad \text{and} \quad
a_2 =2^{-1} \Big(b_1-\sqrt{b_1^2-4 c b_1 b_2^{-1}}\Big)
$$
(reducing, when $c=0$, to $a_1=b_1$ and $a_2=0$, respectively). (Note that the ``if'' part of the assertion follows easily since,
under the relevant conditions, there exist independent 2-component random vectors $\bbb Y^{(1)}$ and $\bbb Y^{(2)}$ so that
$\PP\{\bbb Y^{(1)} = (a_1, a_2 b_1^{-1} b_2)\}=\alpha$, $\PP\{\bbb Y^{(1)}= (a_2, a_1 b_1^{-1} b_2)\}=1-\alpha$,
$\PP\{\bbb Y^{(2)} = (0,0)\}=\beta$, $\PP\{\bbb Y^{(2)} = (-b_1, -b_2)\}=1-\beta$ and $\bbb Y^{(1)} + \bbb Y^{(2)}$ is distributed as $(U,W)$.)
Clearly, the characterization met here implies that the distribution of $(U,W)$ is indecomposable if $U$ and $W$ are nonnegative,
a result referred to in Remark 3 (ii), and also that the distributions of $\bbb X$ appearing in both Examples 4A and 4B are indeed 
indecomposable and hence non-i.d; to see the validity of the claim concerning the distributions of $\bbb X$, note, in particular, that each of these examples has $\EE(|X_1|)=\infty$ with $\bbb X$ meeting the assumptions relative to $(U,W)$.

\medskip
{\bf Remark 4.} For any $p$-component random vector $\bbb X =
(X_1,X_2,\ldots,X_p)$, $p >1$, with d.f. $F$, irrespectively of
whether or not it is i.d., (\ref{eq:propTh2}) is met if and only if
the support of $F$ includes a point $(c_1,c_2,\ldots,c_p)$ so that
$\prod_{r=1}^{p} |c_r| >0$. Suppose now that, for some positive
integer $j$,
$$
\bbb X \stackrel{d}{=} \sum_{n=0}^j \bbb Y^{(n)},
$$
where $\bbb Y^{(n)}$, $n=0,1,\ldots,j$ are independent $p$-component
random vectors such that, for each $n$, given any support point
$(c_1,c_2,\ldots,c_p)$ of the d.f. of $\bbb Y^{(n)}$, there exists a
support point $(d_1,d_2,\ldots,d_p)$ of the d.f. of $\sum_{n' (\neq
n) =0}^j \bbb Y^{(n')}$, for which $\prod_{r=1}^{p} |d_r| \neq 0$
and $c_r d_r>0$ for each $r$ with $c_r \neq 0$. In this case,
$\prod_{r=1}^{p} |c_r+d_r|
> 0$, from which we can see that (\ref{eq:propTh2}) is met, and,
following essentially the relevant steps in the proof of Theorem 2,
appearing (\ref{eq:fanisNewN}) onwards, we can further see that
(\ref{eq:mom1Th2}) implies (\ref{eq:mom1}) with $\bbb Y^{(n)}$ in
place of $\bbb X$ respectively for $n=0,1,\ldots,j$, and hence, in
view of the inequality referred to under brackets below
(\ref{eq:bound1}), in the proof of Theorem 1, that the first
assertion of Theorem 2 holds; this is so irrespectively of whether
or not $\bbb X$ is i.d. (To see, especially, that if
(\ref{eq:mom1Th2}) holds for this $\bbb X$, then, for each $n$,
$\bbb Y^{(n)}$ satisfies (\ref{eq:mom1}), it is sufficient, by
symmetry, taking, in obvious notation, without loss of generality
that $\PP \{ \bbb Y^{(n)} \in (\{0\})^k \times (0,\infty)^{p-k}\} >0$, to
check that this is so for a vector with distribution the same as the
conditional distribution of $\bbb Y^{(n)}$ given that $\bbb Y^{(n)}
\in (\{0\})^k \times (0,\infty)^{p-k}$.) Consequently, we are led
now to a new version of Theorem 2.

\medskip

The following corollary, which is obvious from Theorem 1 and 2,  
presents an extended version of
the first assertion (i.e., the crucial assertion) of Theorem 5.1 of
Gupta {\em et al.} (2009).

\medskip

{\bf Corollary 2.} Let $\bbb X =(X_1,X_2,\ldots,X_p)$ and $\beta_r$,
$r=1,2,\ldots,p$, be as in Theorems 1 and 2. Then,
(\ref{eq:mom1Th2}) is equivalent to (\ref{eq:mom2}).

\medskip

{\bf Remark 5.} In statistical literature, whenever the covariance
between two random variables is mentioned, the relevant random
variables are often assumed to be square-integrable. The following
simple example, which indeed is a slight variation of Example 2 met above, 
illustrates that the covariance may be well-defined
even when this standard assumption does not hold and the joint distribution of the random variables is i.d.:

\medskip

{\bf Example 5.} Let $(X_1,X_2)$ be an i.d. random vector with ch.f.
$\phi$ such that
$$
\phi(\bbb t) = \exp \{-\lambda + \lambda \psi(\bbb t) \}, \quad \bbb
t = (t_1, t_2) \in \RR^2,
$$
where $\lambda \in (0, \infty)$ and $\psi$ is the ch.f. of a random
vector $(V,V^{-1})$ with $V$ as a positive random variable such that
its square equals the modulus of a standard Cauchy random variable.
Then, we have
$$
\EE \left( |X_1X_2| \right) = \EE \left( X_1X_2 \right) < \infty
\quad \text{and} \quad \EE \left( |X_r| \right) = \EE \left( X_r
\right) < \infty, \quad r=1,2.
$$
However, in this case, it is obvious that $ \EE\left( X_1^2
\right)=\EE\left( X_2^2 \right) = \infty$, hence we have the
validity of our claim. (A closer scrutiny of the Example 5 above
tells us, actually, in view of Theorems 1 and 2, something more.
Indeed, the random variables $X_1$ and $X_2$ that we have considered
in this example satisfy the condition that $ \EE\left(
|X_1|^{\alpha_1}|X_2|^{\alpha_2} \right) = \EE\left(
X_1^{\alpha_1}X_2^{\alpha_2} \right) < \infty$ for each $\alpha_1,
\alpha_2 \in (0,2]$.)
Incidentally, it may be worth noting here that, in view of the corollary to Theorem 2 of Shanbhag (1988), implied in Remark 3 (ii) above, any random vector relative to ch.f. $\psi$ of the present example gives us an example of a random vector for which the conclusion of Remark 5 is valid with
``i.d.'' replaced by ``indecomposable''.

\section{\large A characteristic property relative to multivariate self-decomposability}
\setcounter{equation}{0}

Following Urbanik (1969, p. 92), let us define a distribution on
$\RR^p$ ($p \geq 1$) to be self-decomposable (s.d.) if the
corresponding ch.f. $\phi$ satisfies the condition that, for each $c
\in (0,1)$,
\begin{equation}
\label{eq:sd1} \phi(\bbb t) = \phi(c \bbb t) \phi_{c}(\bbb t), \quad
\bbb t \in \RR^p,
\end{equation}
where $\phi_c$ is a ch.f. Essentially, from Sections 8 and 11 of
Chapter XVII of Feller (1966), it is then clear that if
(\ref{eq:sd1}) is met then both $\phi$ and $\phi_c$ are i.d.
ch.f.'s; the partial information that $\phi$ and $\phi_c$ are
nonvanishing also follows, in effect, from the argument given in
Lukacs (1970, p. 162) to show that this is so in the univariate
case. (As a by-product of this, it follows that if $\phi$ is
nonvanishing or, in particular, i.d., then $\phi$ is a s.d. ch.f. if
and only if $\phi(\bbb t) / \phi(c \bbb t)$, $t \in \RR^p$, is an
i.d. ch.f. for each $c \in (0,1)$.) If a distribution on $\RR^p$ is
s.d., we refer to the corresponding ch.f., or a random vector with
this distribution, also as s.d.

\medskip

As pointed out in Section 1, the concept of generalized hyperbolic
distributions is well documented, especially in the univariate case.
Shanbhag \& Sreehari (1979, p. 24) have claimed that there exist
members in the class of multivariate generalized hyperbolic
distributions with $\bbb \bt \neq 0$ (where the distribution
referred to is absolutely continuous with probability density
function (p.d.f.) given by (7.3) of Barndorff-Nielsen (1977)) that
are not s.d. The original example implied in the cited paper to
illustrate this remained unpublished though its version was later
revisited in Pestana (1978, p. 54). More recently, Rao \& Shanbhag
(2004, Example 3.3, Remark 3.4) have produced a simpler example in
support of the claim and have made some further relevant comments on
the issue. There is also Example 6.8 in Gupta {\em et al.} (2006)
that enlightens one with certain related information.

\medskip

Before discussing the results of the present section, we may give
two further pieces of information that are of relevance to these:

(i) Although Feller (1966, XII.8) is concerned with properties of
s.d. distributions in the univariate case, that the properties of
$\phi$ and $\phi_c$ referred to above are met in the multivariate
case is essentially implied by Feller (1966, XII.11); also, it is of
interest to note here that, if $\phi$ satisfies (\ref{eq:sd1}), then
every linear combination of the components of the corresponding
random vector has its ch.f. satisfying the univariate version of
(\ref{eq:sd1}), implying immediately that $\phi$ and $\phi_c$ (of
the multivariate case) are nonvanishing.

\medskip

(ii) To be more precise, in the notation of Section 1, we refer in
this paper to mixtures of $N_n(\mu+u\beta\Delta, u\Delta)$ as
generalized hyperbolic if the mixing distribution relative to $u$ is
generalized inverse Gaussian, i.e., if it has p.d.f. of the form
$$
f(u \mid \lambda, \xi, \psi) = C(\lambda, \xi, \psi) u^{\lambda -1}
\exp \{-(\xi u^{-1} + \psi u)/2 \}, \quad u >0,
$$
with $\xi$, $\psi \geq 0$ for which $\max \{\xi, \psi \} >0$, $C$ as
the normalizing constant depending on the modified Bessel function
of the third kind, and $\lambda$ so that $\int_{0}^{\infty}
u^{\lambda-1} \exp \{-(\xi u^{-1} + \psi u)/2\} du < \infty$.
Obviously, the class of hyperbolic distributions is a subclass of
the class of these distributions, see, e.g., Barndorff-Nielsen
(1977).

\medskip

The following theorem and its corollary subsume some of the major
observations that are made by the examples cited above.

\medskip

{\bf Theorem 3.} Let $V$ be a positive random variable such that
\begin{equation}
\label{eq:sd2} \EE \left(\eee^{s V} \right) = \exp
\left\{\int_{(0,\infty)} \left(\eee^{s v}-1
\right)v^{-(\alpha+1)}g(v)dv \right\}, \quad \text{Re}(s) \leq 0,
\end{equation}
with $\alpha \in [0,1)$ and $g$ as a bounded decreasing nonnegative
real function on $(0, \infty)$ satisfying
$$
\int_{(0,\infty)}\frac{v^{-\alpha}g(v)}{1+v}dv < \infty.
$$
Also, let $\bbb X = (X_1,X_2,\ldots,X_{p-1})$, $p \geq 2$, be a
$(p-1)$-component random vector independent of $V$, with $X_r$,
$r=1,2,\ldots,p-1$, as independent (nondegenerate) symmetric stable
random variables with characteristic exponents $\gamma_r \in [1,2]$,
$r=1,2,\ldots,p-1$, respectively. Then, the $p$-component random
vector
\begin{equation}
\label{eq:sdnn} {\bbb
Z}=\bigg(V,V^{\frac{1}{\gamma_1}}X_1,\ldots,V^{\frac{1}{\gamma_{p-1}}}X_{p-1}\bigg)
\end{equation}
is s.d. if and only if
\begin{equation}
\label{eq:sdn} \alpha -p +1 + \sum_{r=1}^{p-1} \gamma_r^{-1} \geq 0.
\end{equation}

\medskip

{\bf Proof.} In view of (\ref{eq:sd2}), it follows that the ch.f. of
${\bbb Z}$ given by (\ref{eq:sdnn}) is of the form
\begin{equation}
\label{eq:sd3} \phi(\bbb t) = \exp \left\{\int_{(0,\infty)}
\left(\eee^{\iii
t_1v-v\sum_{r=1}^{p-1}\lambda_r|t_{r+1}|^{\gamma_r}}-1
\right)v^{-(\alpha+1)}g(v)dv \right\}, \quad \bbb t
=(t_1,t_2,\ldots,t_{p}) \in \RR^{p},
\end{equation}
with $\lambda_r > 0$ for all $r \in \{1,2,\ldots,p-1\}$; this is
clear since the ch.f. of each symmetric stable random variable
$X_{r}$, with characteristic exponent $\gamma_r$, is of the form
$\phi_{r}(t) = \exp \left\{ -\lambda_r |t|^{\gamma_r}  \right\}$, $t
\in \RR$, and we have
$$
\phi(\bbb t) = \EE \left( \eee^{\iii t_1 V} \prod_{r=1}^{p-1} \EE
\left[ \left( \eee^{\iii t_{r+1} V^{\frac{1}{\gamma_r}}X_r} \right)
\mid V \right]\right), \quad \bbb t \in \RR^{p}.
$$
From (\ref{eq:sd2}), it is obvious that $\phi$ is i.d. Also, on
appealing to Fubini's theorem, it is now clear that the L\'evy
measure in the present case is concentrated on $(0,\infty) \times
\RR^{p-1}$ and is absolutely continuous with respect to Lebesgue
measure on $\RR^{p}$ with Radon-Nikodym derivative $h$ such that the
restriction to $(0,\infty) \times \RR^{p-1}$ of $h$ is given by
\begin{equation}
\label{eq:sd4} h(\bbb y) = y_1^{-(\alpha+1)}g(y_1) \prod_{r=1}^{p-1}
\left[ f_r\left(y_{r+1}/y_1^{\frac{1}{\gamma_r}}\right)
y_1^{-\frac{1}{\gamma_r}} \right], \quad \bbb y =
(y_1,y_2,\ldots,y_{p}) \in (0,\infty) \times \RR^{p-1},
\end{equation}
where $f_r$ denotes the p.d.f. of $X_r$ and is implied by the
inversion theorem to be bounded (continuous). Recalling then that
$\phi$ is s.d. if and only if, for each $c \in (0,1)$, $\phi(\bbb
t)/\phi(c \bbb t)$, $\bbb t \in \RR^p$, is i.d., we can claim that
${\bbb Z}$ is s.d. if and only if
\begin{equation}
\label{eq:sd5} h(\bbb y) \geq c^{-p} h(\bbb y/c), \quad \bbb y \in
(0,\infty) \times \RR^{p-1},
\end{equation}
for all $c \in (0,1)$. Since each symmetric stable distribution is
unimodal with vertex 0 (see, e.g., Lemma 5.10.1 of Lukacs (1970)),
(\ref{eq:sd4}) then implies that (\ref{eq:sd5}) is met for the
required $c$ and $\bbb y$ if and only if (\ref{eq:sdn}) is valid;
observe in particular that if (\ref{eq:sdn}) is not valid, then
(\ref{eq:sd5}) is violated for $\bbb y$ with
$y_{r+1}=o(y_1^{1/\gamma_r})$, $r=1,2,\ldots,p-1$, and $y_1$
sufficiently small.  Hence, Theorem 3 follows. \hfill $\Box$

\medskip

{\bf Remark 6.} In view of (\ref{eq:sd2}) directly, or, trivially,
as a corollary to Theorem 3, it follows that the random variable $V$
considered in Theorem 3 is s.d. Positive stable (with left extremity
zero), gamma and inverse Gaussian random variables provide us with
some specialized versions of $V$ met here; in these cases, we have
$g(v) \propto e^{-\lambda v}$ with $\lambda=0$ for stable and
$\lambda
>0$ otherwise, and, also, have the parameter $\alpha$ respectively
as positive, equal to 0, and  equal to 1/2. Also, in view of the
closure property (under weak convergence) of the class of s.d.
distributions on $\RR^p$, as a corollary to Theorem 3, it now
follows that the ``if'' part of the theorem referred to remains
valid without the assumption that $g$ be bounded.

\medskip

{\bf Corollary 3.} Let $V$ and $X_r$, $r=1,2,\ldots,p-1$, be as in
Theorem 3, and, additionally, suppose that we have $p \geq 3$ and
$\alpha -p + 2 + \sum_{r=1}^{p-2} \gamma_r^{-1} < 0. $ Then, there
exist real $c_r$, $r=1,2,\ldots,p-1$, such that the
$(p-1)$-component random vector ${\bbb W}$ given by
$$
{\bbb W} = \bigg(c_1 V + V^{\frac{1}{\gamma_1}}X_1, \ldots, c_{p-1}
V + V^{\frac{1}{\gamma_{p-1}}}X_{p-1} \bigg)
$$
is not s.d.

\medskip

{\bf Proof.} Clearly, it follows via a standard argument that the
sequence of the random vectors $\{{\bbb W_n}: n=1,2,\ldots\}$,
where, for each $n \geq 1$,
$$
{\bbb W_n} = \bigg(\frac{1}{n}V + V^{\frac{1}{\gamma_1}}X_1, \ldots,
\frac{1}{n} V + V^{\frac{1}{\gamma_{p-2}}}X_{p-2},V+\frac{1}{n}
V^{\frac{1}{\gamma_{p-1}}}X_{p-1}\bigg),
$$
converges in probability and hence in distribution to the random
vector ${\bbb W}^{\star}$ given by
$$
{\bbb W}^{\star} = \bigg(V^{\frac{1}{\gamma_1}}X_1, \ldots,
V^{\frac{1}{\gamma_{p-2}}}X_{p-2},V\bigg).
$$
Since, by Theorem 3, we have that ${\bbb W}^{\star}$ is not s.d.,
appealing to the closure property (under weak convergence) of the
class of s.d. distributions on $\RR^{p-1}$ ($p \geq 3$), we can
readily conclude that, for some $n \geq 1$, ${\bbb W}_n$ is not s.d.
This, in turn, implies that the assertion of the corollary is true;
note that if, for some $n$, $W_n$ is not s.d., then so also is
$(\frac{1}{n} V + V^{\frac{1}{\gamma_{1}}}X_{1},\ldots,\frac{1}{n} V
+V^{\frac{1}{\gamma_{p-2}}}X_{p-2},nV+V^{\frac{1}{\gamma_{p-1}}}X_{p-1})$.
Hence, Corollary 3 follows. \hfill $\Box$

\medskip

{\bf Remark 7.} Although, the problem of finding non-s.d.
multivariate hyperbolic and multivariate generalized hyperbolic
distributions of Barndorff-Nielsen (1977), touched upon in this
article, prompted us to establish Theorem 3, that this latter
theorem is of interest in its own right, is obvious. However, it may
be worth emphasizing here that Corollary 3, which is a corollary to
Theorem 3, identifies certain members of the class of multivariate
generalized hyperbolic distributions of Barndorff-Nielsen (1977),
that are non-s.d.; this follows on noting especially that the
specialized version of the corollary in the case of $\alpha \in
\{0,1/2\}$, $g(v)=\exp\{-\lambda v\}$, $v>0$, and
$\gamma_1=\gamma_2=\ldots\gamma_{p-1}=2$, concerns these
distributions. The examples of non-s.d. distributions given by
Corollary 3 are obviously in the spirit  of those discussed earlier
in Pestana (1978, p. 54) and Rao \& Shanbhag (2004, Example 3.3 \&
Remark 3.4, pp. 2882--2883).

\medskip

{\bf Corollary 4.} Given an integer $p \geq 2$, there exists a
$p$-component random vector of the form of (\ref{eq:sdnn}) that is
not s.d. such that all its lower dimensional subvectors are s.d.

\medskip

{\bf Proof.} Given $p \geq 2$, choose, e.g., $\gamma \in (1, 2]$ and
$\alpha \in [0, 1)$ such that $\alpha = (p-2)(1-\gamma^{-1})$ and,
hence, satisfying also that $\alpha-(p-1)(1-\gamma^{-1}) <0$. Then,
Theorem 3 implies that in the special case of $\gamma_r = \gamma$,
$r=1,2,\ldots,p-1$, with $\alpha$ as stated, the random vector
${\bbb Z}$ is not s.d., but, for each $r \in \{2,3,\ldots,p-1\}$,
the $(p-1)$-component subvector of $\bbb Z$ that does not include
the $r$-th component of $\bbb Z$ is s.d. Also, since the specialized
version in this case of $\bbb Z$ with its first component deleted
satisfies (\ref{eq:sd1}) trivially for each $c \in (0,1)$, it is
obvious that this is s.d. Hence we have the corollary. \hfill $\Box$

\medskip

{\bf Corollary 5.} If $\alpha=0$, the random vector ${\bbb Z}$ is
s.d. if and only if $\gamma_1=\gamma_2=\ldots=\gamma_{p-1}=1$, i.e.,
if and only if $X_1, X_2,\ldots,X_{p-1}$ are Cauchy (up to scale
changes) random variables.

\medskip

{\bf Proof.} The result follows trivially from Theorem 3; also, the
``if'' part of the assertion is immediate on noting that the
concerned ch.f. satisfies (\ref{eq:sd1}) for each $c \in (0,1)$.
 \hfill $\Box$

\medskip

{\bf Remark 8.} Clearly, any random vector ${\bbb Z}^{\star}$ is
i.d. if its ch.f. is of the form
$$
v_{\alpha}^{\star}\bigg(\iii t_1-\sum_{r=1}^{p-1} \lambda_r
|t_{r+1}|^{\gamma_r} \bigg), \quad \bbb t=(t_1,t_2,\ldots,t_p) \in
\RR^p,
$$
with, $p>1$, $\alpha \in [0,1 )$, $\lambda_r > 0$ and $\gamma_r \in
(0,2]$ (i.e., $p$, $\alpha$, $\lambda_r$ and $\gamma_r$ are as in
Theorem 3 but for that $\gamma_r$ is now allowed to lie in $(0,1)$),
and $v_{\alpha}^{\star}$ as the function defined by the left hand
side of (\ref{eq:sd2}); note that ${\bbb Z}^{\star}
\stackrel{d}{=}{\bbb Z}$, where $\bbb Z$ is as in (\ref{eq:sdnn}) if
$\gamma_r \in [1, 2]$, $r=1,2,\ldots,p-1$. Moreover, essentially as
in Example 6.8 of Gupta {\em et al.} (2006), it is now seen that
each projection of ${\bbb Z}^{\star}$ is s.d. if
$\gamma_1=\gamma_2=\ldots=\gamma_{p-1}=2$ and, additionally, $g$ in
(\ref{eq:sd2}) is completely monotone, or, equivalently, by
Bernstein's theorem, denotes the Laplace transform of a measure on
the Borel $\sigma$-field of $\RR_{+}$; this follows since in the
present case also we have the distribution of $V$ to be a member of
the generalized gamma convolution family of Bondesson (1992), on
observing that if $g$ is a Laplace transform as above, then so also
is $v^{-\alpha}g(v)$, $v \in (0,\infty)$. However, by Theorem 3, we
have, in this case, obviously ${\bbb Z}^{\star}$ (and hence ${\bbb
Z}$) to be non-s.d., unless $p = 2$ and $\alpha \in [1/2,1)$.

\medskip

{\bf Remark 9.} The example provided in the proof of Corollary 4
gives us the existence of yet another i.d. distribution on $\RR^p$
($p \geq 2$) that is not s.d., of which each projection is s.d.
This, in conjunction with the information provided by Corollary 4
and Remark 8, compares well with certain findings of L$\acute{\rm
e}$vy (1948) and Shanbhag (1975); the results in the cited
references show us, amongst other things, that there exist
indecomposable distributions on $\RR^p$ ($p \geq 2$) with all
marginals (univariate or otherwise) and projections as i.d. Further
material of relevance to the findings referred to here has appeared
or was cited in, e.g., Davidson (1973), Kendall (1973), Shanbhag
(1974, 1976, 1988), Ostrovskii (1986) and Letac (1992).

\medskip

{\bf Remark 10.} Since the class of distributions on $\RR^p$ that
are s.d. is closed under convolution, from an observation in Remark
3.4 of Rao \& Shanbhag (2004, p. 2883), it is clear that any ch.f.
of the form
$$
v_{\alpha}^{\star}\big(\iii c_1 t_1 + \iii c_2 t_2 - \lambda_1 t_1^2
-\lambda_2 t_2^2 \big), \quad \bbb t=(t_1,t_2) \in \RR^2,
$$
with $v_{\alpha}^{\star}$ as in Remark 8, $(c_1, c_2) \neq (0, 0)$
and $\lambda_1, \lambda_2 > 0$, is s.d., provided that $\alpha \in
[1/2, 1)$ and $g$ is completely monotone; to see this, use the fact
that, in this case, $v^{-(\alpha-1/2)}g(v)$, $v \in (0,\infty)$, is
completely monotone. The cited remark of Rao \& Shanbhag (2004) also
gives some further relevant information on the subject.

\medskip

{\bf Remark 11.} In view of Sato (1999, Example 25.10, pp.
162--164), by (\ref{eq:sd4}), it follows by Fubini's theorem (in the
notation of Theorem 3) that if
$\alpha-p+1+\sum_{r=1}^{p-1}\gamma_r^{-1} <0$, then, for each $\beta
\in (0, p-1-\alpha-\sum_{r=1}^{p-1}\gamma_r^{-1})$ and $\theta \in
(\alpha/(\alpha+\beta),1)$, we have a function $G$ on $(0,\infty)$
such that for all $x \in (0, \infty)$,
\begin{equation}
\label{eq:newFanisGFun} G(x) = \int_{(0,x]}
y_1^{(\alpha+\beta+\sum_{r=1}^{p-1}\gamma_r^{-1})\theta}
\bigg(\prod_{r=1}^{p-1} |y_{r+1}|^{-\theta} \bigg) d \nu(\bbb y) <
\infty,
\end{equation}
where $\nu$ is the L$\acute{\rm e}$vy measure relative to the
distribution of $\bbb Z$ (with $\bbb Z$ as in (\ref{eq:sdnn})). For
each $c \in (0,1)$, in obvious notation, the analogue $G^{(c)}$ of
$G$ with respect to ch.f. $\phi(c \bbb t)$, $\bbb t \in \RR^p$, can
easily be seen to be given by
$c^{(\alpha+\beta-p+1+\sum_{r=1}^{p-1}\gamma_r^{-1})\theta}G(x/c)$,
$x \in (0,\infty)$, with $G$ as in (\ref{eq:newFanisGFun});
obviously, since, for each $c \in (0,1)$ and $x >0$, $G^{(c)}(x) >
G(x)$, it then follows that, for none of $c \in (0,1)$, we have, in
this case, $\phi(\bbb t)/\phi(c \bbb t)$, $\bbb t \in \RR^p$, to be
i.d. Consequently, we have now an alternative argument for proving
the ``only if'' part of Theorem 3.

\medskip

{\bf Remark 12.} Suppose, we define (in the notation in
(\ref{eq:sd2}))
$$
\theta^{\star} = \sup \bigg \{\theta >0: \int_{[1, \infty)}
v^{\theta - \alpha -1} g(v) dv < \infty \bigg \}
$$
Then, by Corollary 1, and the result from Sato (1999, Example 25.10,
pp. 162--164), we have, in view of Fubini's theorem, that, for each
$(\alpha_1,\alpha_2,\ldots,\alpha_p)$, $p \geq 2$, such that
$\alpha_r \in (-1, \gamma_{r-1}^{\star})$, $r=2,3,\ldots,p$, and
$\alpha_1+\sum_{r=2}^{p} \alpha_r \gamma_{r-1}^{-1} \in [0,
\theta^{\star})$,
$$
\EE \bigg( Z_1^{\alpha_1} \prod_{r=2}^{p} |Z_r|^{\alpha_r} \bigg) <
\infty,
$$
where $\bbb Z= (Z_1,Z_2,\ldots,Z_p)$, $p \geq 2$, as in
(\ref{eq:sdnn}) but for a modification that, in this case,
$\gamma_r$, $r=1,2,\ldots,p-1$, $p \geq 2$, are allowed to be less
than 1, and $\gamma_r^{\star}$, $r=1,2,\ldots,p-1$, are so that
$\gamma_r^{\star}=\gamma_r$ if $\gamma_r < 2$ and
$\gamma_r^{\star}=\infty$ if $\gamma_r = 2$. In this remark, as in
Remark 11, we have come across arguments essentially in the spirit
of those met in Section 2 of the paper.

\medskip

{\bf Remark 13.}  Theorem 3 does not hold if the assumption that
$\gamma_r \in [1,2]$, $r=1,2,\ldots,p-1$, is dropped. This is
obvious from the following example:

\medskip

{\bf Example 6.} Let $\{ \phi_n:~n=1,2,\ldots \}$ be a sequence of
i.d. ch.f's on $\RR^3$, such that
$$
\phi_n(\bbb t) = v_{\alpha_{0}}^{\star}\bigg(\iii t_1-t_2^2
-\frac{1}{n} |t_3|^\gamma \bigg), \quad \bbb t=(t_1,t_2,t_3) \in
\RR^3,
$$
with $v_{\alpha}^{\star}$ as in Remark 8, $\alpha_0 \in [0, 1/2)$
and $\gamma \in \big(0, 2/(3-2\alpha_0)\big)$. Clearly, the sequence
$\{\phi_n\}$ converges to a ch.f. on $\RR^3$, which, by Theorem 3,
is not s.d.; note that the limiting random vector $(Y_1,Y_2,Y_3)$ in
this case is so that the ch.f. of $(Y_1,Y_2)$ equals
$v_{{\alpha}_0}^{\star}(\iii t_1-t_2^2)$, $(t_1,t_2) \in \RR^2$,
which is indeed non-s.d. by Theorem 3. Consequently, appealing to
the closure property (under weak convergence) of the class of s.d.
distributions on $\RR^3$, we can claim that there exist, for some
large $n$, ch.f.'s $\phi_n$ that are not s.d. on $\RR^3$ despite the
fact that $\alpha_0-2 + 2^{-1}+ \gamma^{-1}
>0$.

\medskip

{\bf Remark 14.} Applying a result of Zolotarev (1954) (which has
appeared also as Theorem 5.8.4 of Lukacs (1970)), in conjunction
with a standard result stating that each stable distribution is
unimodal (an obvious consequence of a result of Yamazato (1978)
stating that each s.d. distribution is unimodal), one can easily
see, with appropriate scrutiny of L$\acute{\rm e}$vy measures, that
there exist i.d. ch.f.'s on $\RR^2$ of the form
$$
\phi(\bbb t) = v_{\alpha}^{\star}\bigg(\iii t_1- \lambda_1
|t_2|^{\gamma} \bigg), \quad \bbb t=(t_1,t_2) \in \RR^2,
$$
with $v_{\alpha}^{\star}$ as in Remark 8, $\lambda_1 >0$ and $\gamma
\in [1/2,1)$, that are not s.d. for certain $\alpha$ and $g$ (e.g.,
if $\alpha < 1-\gamma$ and $g$ is identically equal to a constant).
This sheds further light on the observation of Remark 13.

\section*{\large Acknowledgements}

We are grateful to Professor Pestana for providing us with a copy of
his Ph.D. Thesis and for sending us some references related to
Section 3 of this paper. Also, we would like to thank the two
referees for their useful comments.

\end{document}